\documentclass[12pt]{article}
\usepackage{amsfonts,amscd,a4}
\newtheorem{theo}{Theorem}
\newtheorem{prop}[theo]{Proposition}
\newtheorem{lemma}[theo]{Lemma}
\newtheorem{coro}[theo]{Corollary}

\newtheorem{defi}[theo]{Definition}

\newcommand{\cA}{{\mathcal A}}
\newcommand{\cB}{{\mathcal B}}
\newcommand{\cC}{{\mathcal C}}

\newcommand{\cF}{{\mathcal F}}

\newcommand{\cK}{{\mathcal K}}
\newcommand{\cL}{{\mathcal L}}

\newcommand{\cP}{{\mathcal P}}

\newcommand{\sC}{{\mathbb C}}

\newcommand{\sN}{{\mathbb N}}

\newcommand{\sR}{{\mathbb R}}

\newcommand{\sT}{{\mathbb T}}

\newcommand{\sZ}{{\mathbb Z}}
\newcommand{\qed}{\rule{1ex}{1ex}}
\begin{document}
\title{The Fredholm index of locally compact band-dominated
operators on $L^p (\sR)$}
\author{Vladimir S. Rabinovich, Steffen Roch}
\date{}
\maketitle
\begin{abstract}
We establish a necessary and sufficient criterion for the
Fredholmness of a general locally compact band-dominated operator
$A$ on $L^p(\sR)$ and solve the long-standing problem of computing 
its Fredholm index in terms of the limit operators of $A$. The 
results are applied to operators of convolution type with almost 
periodic symbol.
\end{abstract}
\section{Introduction} \label{s1}
Throughout this paper, let $1 < p < \infty$, and for each Banach
space $X$, let $L(X)$ stand for the Banach algebra of all bounded
linear operators on $X$, $K(X)$ for the closed ideal of the
compact operators, $B_X$ for the closed unit ball of $X$, and
$X^*$ for the Banach dual space of $X$.

For each function $\varphi \in BUC$, the algebra of the bounded
and uniformly continuous functions on the real line $\sR$, and for
each $t > 0$, set $\varphi_t(x) := \varphi (tx)$ and write
$\varphi I$ for the operator on $L^p(\sR)$ of multiplication by
$\varphi$. An operator $A \in L(L^p(\sR))$ is called {\em
band-dominated} if
\[
\lim_{t \to 0} \|A \varphi_t I - \varphi_t A\| = 0
\]
for each function $\varphi \in BUC$. The set $\cB_p$ of all
band-dominated operators forms a closed subalgebra of $L(L^p(\sR))$.
In this paper we will exclusively deal with band-dominated
operators of the form $I+K$ where $I$ is the identity operator and
$K$ is locally compact (which means that $\varphi A$ and $A
\varphi I$ are compact for each function $\varphi \in BUC$ with
bounded support). We write $\cL_p$ for the set of all locally
compact band-dominated operators on $L^p(\sR)$.

The announced Fredholm criterion and the index formula will be
formulated in terms of limit operators. To introduce this notion,
we will need the shift operators
\[
U_k : L^p(\sR) \to L^p(\sR), \quad (U_k f)(x) := f(x-k)
\]
where $k \in \sZ$. Let $h : \sN \to \sZ$ be a sequence which tends
to infinity in the sense that $|h(n)\to \infty$ as $n \to \infty$.
The operator $A_h \in L(L^p(\sR))$ is called the {\em limit
operator of} $A \in L(L^p(\sR))$ {\em with respect to} $h$ if
\[
\lim_{m \to \infty} \|(U_{-h(m)} A U_{h(m)} - A_h) \varphi I\| = 0
\]
and
\[
\lim_{m \to \infty} \|\varphi (U_{-h(m)} A U_{h(m)} - A_h)\| = 0
\]
for each function $\varphi \in BUC$ with bounded support. The set
of all limit operators of a given operator $A \in L(L^p(\sR))$ is
called the {\em operator spectrum} of $A$ and denoted by
$\sigma_{\! op} (A)$. The operator spectrum splits into two
components $\sigma_+(A) \cup \sigma_-(A)$ which collect the limit
operators of $A$ with respect to sequences $h$ tending to $+
\infty$ and to $- \infty$, respectively.

An operator $A \in L(L^p(\sR))$ is said to be {\em rich} or {\em
to possess a rich operator spectrum} if every
sequence $h$ tending to infinity possesses a subsequence $g$ for
which the limit operator $A_g$ exists. The sets of all rich
operators in $\cB_p$ and $\cL_p$ will be denoted by $\cB_p^\$$ and
$\cL_p^\$$.

Let $\chi_+$ and $\chi_-$ stand for the characteristic functions
of the sets $\sR_+$ and $\sR_-$ of the non-negative and negative
real numbers, respectively. The operators $\chi_+ K \chi_- I$ and $\chi_- K
\chi_+ I$ are compact for each operator $K \in \cL_p$. Indeed, let
$\varepsilon > 0$ be arbitrarily given. Since $K$ is
band-dominated, there is a continuous function $f$ which is 1 on
$[0, \, \infty)$ and 0 on $(- \infty, \, -n_\varepsilon]$ with
sufficiently large $n_\varepsilon$ such that $\|fK - KfI\| <
\varepsilon$. Thus,
\[
\|\chi_+ K \chi_- I - \chi_+ Kf \chi_- I\| =
\|\chi_+ (fK - Kf)  \chi_- I\| < \varepsilon.
\]
The operator $\chi_+ Kf \chi_- I$ is compact since $f \chi_-$ has
a bounded support and $K$ is locally compact. Since further
$\varepsilon$ can be chosen arbitrarily small, the compactness of
$\chi_+ K \chi_- I$ follows. The compactness of $\chi_- K \chi_+
I$ can be checked analogously.

This simple observation implies that, for a Fredholm operator of
the form $A = I + K$ with $K \in \cL_p$, the operators $\chi_+ A
\chi_+ I$ and $\chi_- A \chi_- I$, considered as acting on
$L^p(\sR_+)$ and $L^p(\sR_-)$, are Fredholm operators again. We
call
\[
\mbox{ind}_+ A := \mbox{ind} \, (\chi_+ A \chi_+ I) \quad \mbox{and}
\quad \mbox{ind}_- A := \mbox{ind} \, (\chi_- A \chi_- I)
\]
the {\em plus-} and the {\em minus-index of} $A$. Clearly,
\[
\mbox{ind} \, A = \mbox{ind}_+ A + \mbox{ind}_- A.
\]
Recall in this connection that a bounded linear operator $A$ on 
a Banach space $X$ is said to be {\em Fredholm} if its kernel 
$\mbox{ker} \, A$ and its cokernel $\mbox{coker} \, A := X/\mbox{im} 
\, A$ are linear spaces of finite dimension, and that in this case 
the integer
\[
\mbox{ind} \, A := \mbox{dim} \, \mbox{ker} \, A -  \mbox{dim} \,
\mbox{coker} \, A
\]
is called the {\em Fredholm index} of $A$.

Here is the main result of the present paper.
\begin{theo} \label{t1}
Let $A = I+K$ with $K \in \cL_p^\$$. \\[1mm]
$(a)$ The operator $A$ is Fredholm on $L^p(\sR)$ if and only if
all limit operators of $A$ are invertible and if the norms of
their inverses are uniformly bounded. \\[1mm]
$(b)$ If $A$ is Fredholm, then for arbitrary operators $B_+ \in
\sigma_+ (A)$ and $B_- \in \sigma_- (A)$,
\begin{equation} \label{e2}
\mbox{\rm ind}_+ B_+ = \mbox{\rm ind}_+ A \quad \mbox{and} \quad
\mbox{\rm ind}_- B_- = \mbox{\rm ind}_- A
\end{equation}
and, consequently,
\begin{equation} \label{e3}
\mbox{\rm ind} \, A = \mbox{\rm ind}_+ B_+ + \mbox{\rm ind}_- B_-.
\end{equation}
\end{theo}
This result has a series of predecessors. One of the simplest
classes of band-dominated and locally compact operators on
$L^p(\sR)$ is constituted by the operators of convolution by
$L^1(\sR)$-functions and by the restrictions of these operators to
the half line, the classical Wiener-Hopf operators. The theory of
the convolution type operators on the half line originates from
the fundamental papers by Krein and Gohberg/Krein \cite{Kre1,GKr1}
where the Fredholm theory for these operators is established and
an index formula is derived. See also the monograph \cite{GoF1} by
Gohberg/Feldman for an axiomatic approach to this circle of
questions. For convolution type operators with variable
coefficients which stabilize at infinite, a Fredholm criterion and
an index formula have been obtained by Karapetiants/Samko in
\cite{KaS1}; see also their monograph \cite{KaS2}.

In \cite{RRS1,RRS2}, there is developed the limit operator
approach to study Fredholm properties of general band-dominated
operators on spaces $l^p$ of vector-valued sequences. In
\cite{RaR5} we demonstrated that this approach also applies to
operators of convolution type acting on $L^p$ spaces if a suitable
discretization reducing $L^p$- to $l^p$-spaces is performed. (To
be precisely: If the sequences in $l^p$ take their values in an
infinite dimensional Banach space, then we derived in \cite{RRS2}
a criterion for a generalized form of Fredholmness, called
$\cP$-Fredholmness; see below. But the results of \cite{RaR5}
refer to common Fredholmness.) The long standing problem to
determine the Fredholm index of a band-dominated operator in terms
of its limit operators, too, has been finally solved in
\cite{RRR1} for band-dominated operators on the space $l^2$ with
scalar-valued sequences. All mentioned results can be also found
in the monograph \cite{RRS4}. The index formula has been
generalized to $l^p$-spaces in \cite{Roc9}. In the present paper
we will undertake a further generalization to band-dominated
operators with compact entries acting on $l^p$-spaces of
vector-valued functions. Thereby these results will get the right
form to become applicable to locally compact band-dominated
operators on $L^p$-spaces (and thus, to prove assertion $(b)$ of
the theorem).

The paper is organized as follows. We start with recalling some
basic facts on sequences of compact operators. For the reader's
convenience, the proofs are included. The main work will be done
in Section \ref{s3} where we will derive the Fredholm criterion
and the index formula for band-dominated operators on $l^p$ with
compact entries. In Section \ref{s4}, these results will be applied
to locally compact band-dominated operators on $L^p$ which mainly
requires to construct a suitable discretization mapping. Some
applications will be discussed in the final section.

This work was supported by the CONACYT project 43432. The authors 
are grateful for this support.
\section{Sequences of compact operators} \label{s2}
Let $X$ be a complex Banach space which enjoys the following {\em
symmetric approximation property} ($sap$): There is a sequence
$(\Pi_N)_{N \ge 1}$ of projections (= idempotents) $\Pi_N \in
L(X)$ of finite rank such that $\Pi_N \to I$ and $\Pi_N^* \to I^*$
strongly as $N \to \infty$. Evident examples of Banach spaces with
$sap$ are the separable Hilbert spaces, the spaces $l^p(\sZ^K)$
and the spaces $L^p[a, \, b]$. It is also clear that if $X$ is a
reflexive Banach space with $sap$, then $X^*$ has $sap$, too, and
the corresponding projections can be chosen as $\Pi_N^*$.
\begin{defi} \label{d4}
A sequence $(K_n)$ of operators in $L(X)$ is said to be \\[1mm]
$(a)$ {\em relatively compact} if the norm closure of $\{ K_n :
n \in \sN \}$ is compact in $L(X)$; \\[1mm]
$(b)$ {\em collectively compact} if the set $\cup_{n \in \sN}
K_n B_X$ is relatively compact in $X$; \\[1mm]
$(c)$ {\em uniformly left (right, two-sided) approximable} if, for
each $\varepsilon > 0$ there is an $N_0$ such that, for each $n
\in \sN$ and each $N \ge N_0$,
\[
\|K_n - \Pi_N K_n\| < \varepsilon \quad (\|K_n - K_n \Pi_N \| <
\varepsilon, \quad \|K_n - \Pi_N K_n \Pi_N \| < \varepsilon).
\]
\end{defi}
Note that the uniform left approximability of $(K_n)$ is
equivalent to
\[
\lim_{N \to \infty} \sup_{n \in \sN} \|K_n - \Pi_N K_n\| = 0.
\]
\begin{prop} \label{p5}
Let $X$ be a Banach space with $sap$. The following conditions are
equivalent for a sequence $(K_n)$ of compact operators on $X \!:$ \\[1mm]
$(a)$ $(K_n)$ is relatively compact; \\[1mm]
$(b)$ $(K_n)$ and $(K_n^*)$ are collectively compact; \\[1mm]
$(c)$ $(K_n)$ is uniformly left and uniformly right approximable; \\[1mm]
$(d)$ $(K_n)$ is uniformly two-sided approximable.
\end{prop}
{\bf Proof.} $(a) \Rightarrow (b)$: Let $(x_n)$ be a sequence in
$\cup_n K_n B_X$. For each $n \in \sN$, choose $r(n) \in \sN$ and
$y_n \in B_X$ such that $x_n = K_{r(n)} y_n$. By hypothesis $(a)$,
the sequence $(K_{r(n)})$ has a convergent subsequence
$(K_{r(n_k)})$. Let $K$ denote the limit of that subsequence. Then
\begin{equation} \label{e6}
\|x_{n_k} - K y_{n_k} \| = \|K_{r(n_k)} y_{n_k} - K y_{n_k}\| \le
\|K_{r(n_k)} - K\| \to 0.
\end{equation}
Since $K$ is compact and $\|y_{n_k}\| \le 1$, the sequence
$(K_{r(n_k)})$ has a convergent subsequence. From (\ref{e6}) we
conclude that then the sequence $(x_{n_k})$ (hence, the sequence
$(x_n)$) has a convergent subsequence, too. This yields the
collective compactness of the sequence $(K_n)$. Since $(K_n^*)$ is
relatively compact whenever $(K_n)$ is relatively compact, the
collective compactness of $(K_n^*)$ follows in the same way.
\\[1mm]
$(b) \Rightarrow (c)$: We will show that the collective
compactness of $(K_n)$ implies the uniform left approximability of
that sequence. We will not make use of the strong convergence of
$\Pi_N$ to $I^*$ in this part of the proof. So is becomes evident
that then also the collective compactness of $(K_n^*)$ implies the
uniform left approximability of $(K_n^*)$ with respect to the
sequence $(\Pi_N^*)$ which is equivalent to the uniform right
approximability of $(K_n)$.

Contrary to what we want to show, assume that $(K_n)$ is not
uniformly left approximable. Then there are an $\varepsilon > 0$,
a monotonically increasing sequence $(N(r))_{r \ge 1}$ and
operators $K_{n(r)} \in \{ K_n : n \in \sN \}$ such that
\[
\|(I - \Pi_{N(r)}) K_{n(r)} \| \ge \varepsilon \quad \mbox{for
all} \; r \in \sN.
\]
Choose $x_{n(r)} \in B_X$ such that
\begin{equation} \label{e7}
\|(I - \Pi_{N(r)}) K_{n(r)} x_{n(r)}\| \ge \varepsilon/2 \quad
\mbox{for all} \; r \in \sN.
\end{equation}
By hypothesis $(b)$, the sequence $(K_{n(r)} x_{n(r)})$ has a
convergent subsequence. Let $x_0$ denote its limit. We
conclude from (\ref{e7}) that $\|(I - \Pi_{N(r)}) x_0\| \ge
\varepsilon/4$ for all sufficiently large $r$. Letting $r$ go to
infinity, we arrive at a contradiction. \\[1mm]
$(c) \Rightarrow (d)$: This implication follows immediately from
\begin{eqnarray*}
\|K_n - \Pi_N K_n \Pi_N\| & \le & \|K_n - \Pi_N K_n\| + \|\Pi_N
K_n - \Pi_N K_n \Pi_N\| \\
& \le & \|K_n - \Pi_N K_n\| + \|\Pi_N\| \, \|K_n - K_n \Pi_N\|
\end{eqnarray*}
and from the uniform boundedness of the projections $\Pi_N$ due to
the Banach-Steinhaus theorem. \\[1mm]
$(d) \Rightarrow (a)$: We consider a subsequence of $(K_n)$ which
we write as $(K_n)_{n \in \sN_0}$ with an infinite subset $\sN_0$
of $\sN$. Since the projections $\Pi_N$ have finite rank, there
are an infinite subset $\sN_1$ of $\sN_0$ such that the sequence
$(\Pi_1 K_n \Pi_1)_{n \in \sN_1}$ converges, an infinite subset
$\sN_2$ of $\sN_1$ such that the sequence $(\Pi_2 K_n \Pi_2)_{n
\in \sN_2}$ converges, etc. Thus, for each $N \ge 1$, one finds an
infinite subset $\sN_N$ of $\sN_{N-1}$ such that the sequence
$(\Pi_N K_n \Pi_N)_{n \in \sN_N}$ converges. Let $k(n)$ denote the
$n$th number in $\sN_N$ (ordered with respect to the relation $<$)
and set $\hat{K}_n := K_{k(n)}$. Clearly, $(\hat{K}_n)_{n \ge 1}$
is a subsequence of each of the sequences $(K_n)_{n \in \sN_N}$ up
to finitely many entries. Thus, for each $N \in \sN$, the sequence
$(\Pi_N \hat{K}_n \Pi_N)_{n \ge 1}$ converges. Now we have
\[
\hat{K}_n - \hat{K}_m = (\hat{K}_n - \Pi_N \hat{K}_n \Pi_N) -
(\hat{K}_m - \Pi_N \hat{K}_m \Pi_N) + \Pi_N (\hat{K}_n -
\hat{K}_m) \Pi_N.
\]
Let $\varepsilon > 0$. By hypothesis $(d)$, there is an $N$ such
that
\[
\| \hat{K}_n - \Pi_N \hat{K}_n \Pi_N \| < \varepsilon/3
\]
for all $n \in \sN$. Fix this $N$, and choose $n_0$ such that
\[
\| \Pi_N (\hat{K}_n - \hat{K}_m) \Pi_N \| < \varepsilon/3
\]
for all $m, \, n \ge n_0$ which is possible due to the convergence
of the sequence $(\Pi_N \hat{K}_n \Pi_N)_{n \ge 1}$. Hence,
$\|\hat{K}_n - \hat{K}_m\| < \varepsilon$ for all $m, \, n \ge
n_0$. This implies the convergence of the sequence $(\hat{K}_n)$
and, thus, the relative compactness of $(K_n)$. \hfill \qed
\section{The Fredholm index of discrete band-domi\-nated
operators with compact entries} \label{s3}
Let $X$ be a complex Banach space with $sap$. By $E := l^p(\sZ, \,
X)$ we denote the Banach space of all sequences $x : \sZ \to X$
with
\[
\|x\|_E^p := \sum_{n \in \sZ} \|x_n\|_X^p < \infty.
\]
For $k \in \sZ$, let $V_k : l^p(\sZ, \, X) \to l^p(\sZ, \, X)$
stand for the shift operator $(V_k x)_n := x_{n-k}$. In what follows, 
we will have to consider shift operators on different spaces 
$l^p(\sZ, \, X)$. In order to indicate the underlying space we will 
sometimes also write $V_{k, \, X}$ for the shift operator $V_k$ on 
$l^p(\sZ, \, X)$. Further, for each non-negative integer $n$, let the
projection operators $P_n : l^p(\sZ, \, X) \to l^p(\sZ, \, X)$ be
defined by
\[
(P_n x)_k := \left\{
\begin{array}{lll}
x_k & \mbox{if} & |k| \le n \\
0 & \mbox{if} & |k| > n,
\end{array}
\right.
\]
and set $Q_n := I - P_n$ and $\cP := (P_n)_{n \ge 0}$. Sometimes we
will also write $P_{n, \, X}$ in place of $P_n$ in order to indicate
the underlying space.

Each operator $A \in L(E)$ can be represented in the obvious way
by a two-sided infinite matrix with entries in $L(X)$ (in analogy
with the representation of an operator on $l^p(\sZ) := l^p(\sZ, \,
\sC)$ with respect to the standard basis). The operator $A \in
L(E)$ is called a {\em band operator} if its matrix representation
$(A_{ij})$ is a band matrix, i.e., if there is a $k \in \sN$ such
that $A_{ij} = 0$ if $|i-j| \ge k$. The closure of the set of all
band operators on $E$ is a closed subalgebra of $L(E)$ which we
denote by $\cA_E$. The elements of $\cA_E$ will be called {\em
band-dominated operators}. By $\cC_E$ we denote the closed ideal
of $\cA_E$ which consists of all band-dominated operators which
have only compact entries in their matrix representation.

Following the terminology introduced in \cite{RRS4}, an operator
$K \in L(E)$ is called {\em $\cP$-compact} if
\[
\lim_{n \to \infty} \|K Q_n\| = \lim_{n \to \infty} \|Q_n K\| = 0.
\]
We denote the set of all $\cP$-compact operators by $K(E, \, \cP)$,
and we write $L(E, \, \cP)$ for the set of all operators $A \in
L(E)$ for which both $AK$ and $KA$ are $\cP$-compact whenever $K$
is $\cP$-compact. Then $L(E, \, \cP)$ is a closed subalgebra of
$L(E)$ which contains $K(E, \, \cP)$ as a closed ideal.
\begin{defi} \label{d8}
An operator $A \in L(E, \, \cP)$ is called {\em $\cP$-Fredholm} if
the coset $A + K(E, \, \cP)$ is invertible in the quotient algebra
$L(E, \, \cP)/K(E, \, \cP)$, i.e., if there exist an operator $B
\in L(E, \, \cP)$ and operators $K, \, L \in K(E, \, \cP)$ such
that $B A = I + K$ and $A B = I + L$.
\end{defi}
This definition is equivalent to the following one: An operator $A
\in L(E, \, \cP)$ is $\cP$-Fredholm if and only if there exist an
$m \in \sN$ and operators $L_m, \, R_m \in L(E, \, \cP)$ such that
\[
L_m A Q_m = Q_m A R_m = Q_m.
\]
Thus, $\cP$-Fredholmness is often referred to as {\em local
invertibility at infinity}. If $X$ has finite dimension, then the
notions $\cP$-Fredholmness and Fredholmness are synonymous.

All band-dominated operators belong to $L(E, \, \cP)$. This can be
easily checked for the two basic types of band-dominated
operators: the shift operators and the operators of multiplication
by a function in $l^\infty(\sZ, \, L(X)$), and it follows for
general band-dominated operators since $L(E, \, \cP)$ is a closed
algebra. Hence, it makes sense to speak about their
$\cP$-Fredholmness. A criterion for the $\cP$-Fredholmness of a
band-dominated operator $A$ can be given in terms of the limit
operators of $A$. These are, in analogy with the notions from
Section \ref{s1}, defined as follows. Let $A \in L(E)$, and let $h
: \sN \to \sZ$ be a sequence which tends to infinity. An operator
$A_h \in L(E)$ is called a {\em limit operator of} $A$ {\em with
respect to the sequence} $h$ if
\[
\lim_{n \to \infty} \|P_k (V_{-h(n)} A V_{h(n)} - A_h)\| = \lim_{n
\to \infty} \|(V_{-h(n)} A V_{h(n)} - A_h) P_k\| = 0
\]
for every $k \in \sN$. The set of all limit operators of $A$ will
be denoted by $\sigma_{\! op}(A)$ and is called the {\em operator
spectrum} of $A$ again. An operator $A \in L(E)$ is said to be
{\em rich} or to possess a {\em rich operator spectrum} if each
sequence $h$ which tends to infinity possesses a subsequence $g$
for which the limit operator $A_g$ exists. We refer to the rich
operators in $\cA_E$ as {\em rich} band-dominated operators and
write $\cA_E^\$$ and $\cC_E^\$$ for the Banach algebra of the rich
band-dominated operators and for its closed ideal consisting of
the rich operators in $\cC_E$.

The following is the main result on $\cP$-Fredholmness of rich
band-dominated operators. Its proof is in \cite{RRS4}, Theorem
2.2.1.
\begin{theo} \label{t9}
An operator $A \in \cA_E^\$$ is $\cP$-Fredholm if and only if each
of its limit operators is invertible and if the norms of their
inverses are uniformly bounded, i.e.,
\[
\sup \{ \| (A_h)^{-1} \| : A_h \in \sigma_{\! op}(A) \} < \infty.
\]
\end{theo}
In case $X = \sC$, $\cP$-Fredholmness coincides with common
Fredholmness. In this case one can also express the Fredholm index
of a Fredholm band-dominated operator in terms of the (local)
indices of its limit operators. To cite these results from
\cite{RRR1,Roc9}, let $P : l^p(\sZ, \, X) \to l^p(\sZ, \, X)$
refer to the projection operator
\[
(P x)_k := \left\{
\begin{array}{lll}
x_k & \mbox{if} & k \ge 0 \\
0 & \mbox{if} & k < 0,
\end{array}
\right.
\]
and set $Q := I-P$. If necessary, we will write also $P_X$ in
place of $P$. Then, for each band-dominated operator on $l^p(\sZ,
\, \sC)$, the operators $PAQ$ and $QAP$ are compact. This is
obvious for band operators in which case $PAQ$ and $QAP$ are of
finite rank, and it follows for general band-dominated operators
by an obvious approximation argument. Consequently, the operators
$A - (PAP+Q)(P+QAQ)$ and $A - (P+QAQ)(PAP+Q)$ are compact, which
implies that a band-dominated operator on $l^p(\sZ, \, \sC)$ is
Fredholm if and only if both operators $PAP+Q$ and $P+QAQ$ are
Fredholm and that
\[
\mbox{ind} \, A = \mbox{ind} \, (PAP+Q) + \mbox{ind} \, (P+QAQ).
\]
In this case we call
\[
\mbox{ind}_+ A := \mbox{ind} \, (PAP+Q) \quad \mbox{and} \quad
\mbox{ind}_- A := \mbox{ind} \, (P+QAQ)
\]
the {\em plus-} and the {\em minus-index of} $A$. Finally, let
$\sigma_{\! op} (A) = \sigma_+(A) \cup \sigma_-(A)$, the latter
components collecting the limit operators of $A$ with respect to
sequences $h$ tending to $+ \infty$ and to $- \infty$,
respectively, and note that in case $X=\sC$ all band-dominated
operators are rich.
\begin{theo} \label{t10}
Let $X = \sC$, and let $A$ be a Fredholm band-dominated operator
on $l^p(\sZ)$. Then, for arbitrary operators $B_+ \in \sigma_+
(A)$ and $B_- \in \sigma_- (A)$,
\begin{equation} \label{e11}
\mbox{\rm ind}_+ B_+ = \mbox{\rm ind}_+ A \quad \mbox{and} \quad
\mbox{\rm ind}_- B_- = \mbox{\rm ind}_- A
\end{equation}
and, consequently,
\begin{equation} \label{e12}
\mbox{\rm ind} \, A = \mbox{\rm ind}_+ B_+ + \mbox{\rm ind}_- B_-.
\end{equation}
\end{theo}
In particular, all operators in $\sigma_+ (A)$ have the same
plus-index, and all operators in $\sigma_- (A)$ have the same
minus-index.

It is the goal of the present section to generalize the assertion
of Theorem \ref{t10} to operators acting on $E = l^p(\sZ, \, X)$
with a general Banach space $X$ with $sap$ which are of the form
$I+K$ with $K \in \cC_E^\$$. A first observation is that for these
operators $\cP$-Fredholmness and common Fredholmness coincide.
\begin{prop} \label{p13}
An operator in $I + \cC_E$ is Fredholm if and only if it is
$\cP$-Fredholm.
\end{prop}
{\bf Proof.} We claim that
\begin{equation} \label{e14}
\cC_E \cap K(E, \, \cP) = K(E).
\end{equation}
The inclusion $K(E) \subseteq \cC_E$ is evident, and the inclusion
$K(E) \subseteq K(E, \, \cP)$ holds since the projections $P_n$
and $P_n^*$ converge strongly to the identity operators on $E$ and
$E^*$, respectively. Thus, $K(E) \subseteq \cC_E \cap K(E, \,
\cP)$. For the reverse inclusion, let $K \in \cC_E \cap K(E, \,
\cP)$. Since $K \in \cC_E$, one has $P_nK \in K(E)$ for every $n$,
and since $K \in K(E, \, \cP)$, one has $\|K - P_nK\| \to 0$.
Thus, $K \in K(E)$, which verifies (\ref{e14}).

Since $K(E) \subseteq K(E, \, \cP)$ by (\ref{e14}), every
Fredholm operator in $L(E, \, \cP)$ is $\cP$-Fredholm.
For the reverse implication, let $A := I + K$ with $K \in \cC_E$
be a $\cP$-Fredholm operator. Then there are operators $B \in L(E,
\, \cP)$ and $L \in K(E, \, \cP)$ such that $BA = I - L$. Set $R
:= I - KB$. Then
\[
RA - I = A - I - KBA = K - KBA = K(I - BA) = KL.
\]
Since $KL \in \cC_E \cap K(E, \, \cP)$ is compact by (\ref{e14}),
the operator $R$ is a left Fredholm regularizer for $A$. Similarly
one checks that $A$ possesses a right Fredholm regularizer. Thus,
the operator $A$ is Fredholm. \hfill \qed \\[3mm]
Combining Proposition \ref{p13} with Theorem \ref{t9} one gets the
following.
\begin{coro} \label{c15}
Let $A := I + K$ with $K \in \cC_E^\$$. Then the operator $A$ is
Fredholm if and only if each of its limit operators is invertible
and if the norms of their inverses are uniformly bounded.
\end{coro}
We will make use of the following lemma several times.
\begin{lemma} \label{l16}
Every band-dominated operator in $\cC_E$ (resp. in $\cC_E^\$$) is
the norm limit of a sequence of band operators $\cC_E$ (resp. in
$\cC_E^\$$).
\end{lemma}
This can be proved in exactly the way as we derived Theorem 2.1.18
in \cite{RRS4} which states that every rich band-dominated
operators is the norm limit of a sequence of rich band operators.
\hfill \qed \\[3mm]
As a first consequence of the $\cC_E$-version of Lemma \ref{l16}
we conclude that $PAQ$ and $QAP$ are compact operators for each
operator $A \in I + \cC_E$. Indeed, this is obvious for $A$ being
a band operators in which case $PAQ$ and $QAP$ have only a finite
number of non-vanishing entries, and these are compact. The case
of general $A$ follows by an obvious approximation argument.
Consequently, the operators $A - (PAP+Q)(P+QAQ)$ and $A -
(P+QAQ)(PAP+Q)$ are compact, which implies that an operator $A \in
I + \cC_E$ is Fredholm if and only if both operators $PAP+Q$ and
$P+QAQ$ are Fredholm. In this case, the integers
\[
\mbox{ind}_+ A := \mbox{ind} \, (PAP+Q) \quad \mbox{and} \quad
\mbox{ind}_- A := \mbox{ind} \, (P+QAQ)
\]
are called the {\em plus-} and the {\em minus-index of} $A$.
Clearly,
\begin{equation} \label{e17}
\mbox{ind} \, A = \mbox{ind}_+ A + \mbox{ind}_- A.
\end{equation}
Finally, let $\sigma_{\! op} (A) = \sigma_+(A) \cup \sigma_-(A)$
in analogy with the case $X = \sC$.

Here is the announced result for the indices of Fredholm operators
in $I + \cC_E^\$$.
\begin{theo} \label{t18}
Let $A \in I + \cC_E^\$$ be a Fredholm operator. Then, for
arbitrary operators $B_+ \in \sigma_+ (A)$ and $B_- \in \sigma_-
(A)$,
\begin{equation} \label{e19}
\mbox{\rm ind}_+ B_+ = \mbox{\rm ind}_+ A \quad \mbox{and} \quad
\mbox{\rm ind}_- B_- = \mbox{\rm ind}_- A
\end{equation}
and, consequently,
\begin{equation} \label{e20}
\mbox{\rm ind} \, A = \mbox{\rm ind}_+ B_+ + \mbox{\rm ind}_- B_-.
\end{equation}
\end{theo}
The remainder of this section is devoted to the proof of Theorem
\ref{t18}. We will verify this theorem by reducing its assertion
step by step until we will arrive at operators on $l^p(\sZ, \,
\sC)$ (with scalar entries) for which the result is known (Theorem
\ref{t10}). The first step of the reduction procedure is based on
the following observation.
\begin{prop} \label{p21}
Let $\cF$ be a dense subset of the set of all Fredholm operators
in $I + \cC_E^\$$. If the assertion of Theorem $\ref{t18}$ holds
for all operators in $\cF$, then it holds for all Fredholm
operators in $I + \cC_E^\$$.
\end{prop}
{\bf Proof.} Let $A$ be a Fredholm operator in $I + \cC_E^\$$, and
let $B \in \sigma_+ (A)$. We will show that
\begin{equation} \label{e22}
\mbox{ind}_+ B = \mbox{ind}_+ A,
\end{equation}
which settles the plus-assertion of (\ref{e19}). The
minus-assertion follows similarly, and (\ref{e19}) implies
(\ref{e20}) via (\ref{e17}).

To prove (\ref{e22}), choose a sequence $(A_n)$ of operators in
$\cF$ which converges to $A$ in the operator norm, and let $h$ be
a sequence tending to $+ \infty$ such that $B = A_h$. Employing
Cantor's diagonal method, we construct a subsequence $g$ of $h$
for which all limit operators $(A_n)_g$ exist. For the details of
this construction, consult the proof of Proposition 1.2.6 in
\cite{RRS4}. From Proposition 1.2.2 $(e)$ in \cite{RRS4} we
conclude that $\|B - (A_n)_g\| = \|A_g - (A_n)_g\| \to 0$. Now one
has
\[
\mbox{ind}_+ (A_n)_g = \mbox{ind}_+ A_n \quad \mbox{for all} \; n
\in \sN
\]
and this implies (\ref{e22}) by letting $n$ go to infinity due to
the continuity of the index. \hfill \qed \\[3mm]
Our choice of the set $\cF$ is as follows. The $\cC_E^\$$-version
of Lemma \ref{l16} allows one to approximate each band-dominated
operator $A = I + K$ with $K \in \cC_E^\$$ by a sequence of band
operators $A_n := I + K_n$ with $K_n \in \cC_E^\$$. Each band
operator $K_n$ can be written as a sum
\begin{equation} \label{e23}
K_n = \sum_{k \in \sZ} K_n^{(k)} V_k
\end{equation}
with only finitely many non-vanishing items. The coefficients
$K_n^{(k)}$ in (\ref{e23}) are operators of multiplication by
sequences of compact operators on $X$, and these multiplication
operators are rich whenever $K_n$ is rich. From Theorem 2.1.16 in
\cite{RRS4} we know that a multiplication operator is rich if and
only if the set of its entries is relatively compact in $L(X)$. So
we conclude from the equivalence between $(a)$ and $(d)$ in
Proposition \ref{p5} that each coefficient $K_n^{(k)}$ in
(\ref{e23}) can be approximated as closely as desired by a
sequence $(K_{n, \, N}^{(k)})_{N \in \sN}$ of multiplication
operators the entries of which map $\mbox{im} \, \Pi_N$ into
itself and act on $\mbox{im} \, (I_X - \Pi_N)$ as the zero
operator. Thus, one can approximate the operator $A = I + K$ as
closely as desired by band operators $A_{n, \, N} = I + K_{n, \,
N}$ where the entries of $K_{n, \, N}$ map $\mbox{im} \, \Pi_N$
into itself and act as the zero operator on $\mbox{im} \, (I_X -
\Pi_N)$. We denote the set of all operators $K_{n, \, N}$ of this
form by $\cC_{E, \, N}$. Note that the operators in $\cC_{E, \,
N}$ are automatically rich.

Further, if $A = I + K$ is a Fredholm operator, then the operators
$A_{n, \, N} = I + K_{n, \, N}$ are Fredholm for all sufficiently
large $n$ and $N$. Thus, we can choose $\cF$ as the set of all
Fredholm operators $I + K_{n, \, N}$ with $K_{n, \, N} \in \cC_{E,
\, N}$. By Proposition \ref{p21}, it remains to prove Theorem
\ref{t18} for these operators.

We agree upon writing $X_N$ in place of $\mbox{im} \, \Pi_N$ if we
want to consider $\mbox{im} \, \Pi_N$ as a Banach space in its own
right, not as a subspace of $X$. Further we introduce the mappings
\[
R : l^p(\sZ, \, X) \to l^p(\sZ, \, X_N), \quad (x_n) \mapsto
(\Pi_N x_n)
\]
where $\Pi_N x_n$ is considered as an element of $X_N$, and
\[
L : l^p(\sZ, \, X_N) \to l^p(\sZ, \, X), \quad (x_n) \mapsto (x_n)
\]
where the $x_n$ on the right-hand side are considered as elements
of $X$. Clearly, $RL$ is the identity operator on $l^p(\sZ, \,
X_N)$, whereas $LR$ is the projection
\[
\Pi : l^p(\sZ, \, X) \to l^p(\sZ, \, X), \quad (x_n) \mapsto
(\Pi_N x_n),
\]
now with the $\Pi_N x_n$ being considered as elements of $X$. We
are going to show that the operators $A = I + K_{n, \, N}$ as well
as their limit operators behave well under the mapping $A \mapsto
RAL$.
\begin{prop} \label{p24}
Let $A = I + K_{n, \, N}$ with $K_{n, \, N} \in \cC_{E, \, N}$.
\\[1mm]
$(a)$ If $A$ is a Fredholm operator on $l^p(\sZ, \, X)$, then
$RAL$ is a Fredholm operator on $l^p(\sZ, \, X_N)$, and the
Fredholm indices of $A$ and $RAL$ coincide. \\[1mm]
$(b)$ If the limit operator of $A$ with respect to a sequence $h :
\sN \to \sZ$ exists, then the limit operator of $RAL$ with respect
to $h$ exists, too, and $(RAL)_h = RA_hL$.
\end{prop}
{\bf Proof.} Since $A$ is Fredholm, there are operators $B, \, T$
on $l^p(\sZ, \, X)$ with $T$ compact such that
\begin{equation} \label{e25}
BA = I + T.
\end{equation}
For $x \in \mbox{ker} \, A$ one gets $x + Tx = 0$, whence $x \in
\mbox{im} \, T$. Hence, $\mbox{dim} \, \mbox{ker} \, A \le
\mbox{rank} \, T$ for each pair $(B, \, T)$ such that (\ref{e25})
holds. One can choose the pair $(B, \, T)$ even in such a way that
$\mbox{dim} \, \mbox{ker} \, A = \mbox{rank} \, T$. For write $X$
as a direct sum $\mbox{ker} \, A \oplus X_0$ and let
$P_{\mbox{\scriptsize ker} \, A}$ refer to the projection from $X$
onto $\mbox{ker} \, A$ parallel to $X_0$. Then
\[
A (I - P_{\mbox{\scriptsize ker} \, A}) : \mbox{im} \, (I -
P_{\mbox{\scriptsize ker} \, A}) \to \mbox{im} \, A
\]
is an invertible operator. Let $B$ denote its inverse. Then $BA (I
- P_{\mbox{\scriptsize ker} \, A}) = I - P_{\mbox{\scriptsize ker}
\, A}$ and
\[
BA = I - P_{\mbox{\scriptsize ker} \, A} + BA P_{\mbox{\scriptsize
ker} \, A} = I - (I-BA) P_{\mbox{\scriptsize ker} \, A}.
\]
Clearly, $\mbox{rank} \, (I-BA) P_{\mbox{\scriptsize ker} \, A}
\le \mbox{dim} \, \mbox{ker} \, A$. Thus, one can indeed assume
that (\ref{e25}) holds with $\mbox{dim} \, \mbox{ker} \, A =
\mbox{rank} \, T$. From (\ref{e25}) we get
\[
RBAL = RL + RTL = I + RTL,
\]
and since $L = \Pi L$ and $A$ commutes with $\Pi$, we obtain
\begin{equation} \label{e26}
RBL \, RAL = I + RTL.
\end{equation}
In the same way, $AB = I + T^\prime$ with $T^\prime$ compact
implies that $RAL \, RBL = I + RT^\prime L$ with $RT^\prime L$
compact. Hence, $RAL$ is Fredholm, and (\ref{e26}) moreover shows
that
\[
\mbox{dim} \, \mbox{ker} \, RAL \le \mbox{rank} \, RTL \le
\mbox{rank} \, T = \mbox{dim} \, \mbox{ker} \, A.
\]
For the reverse estimate, let $B, \, T$ be operators on $l^p(\sZ,
\, X_N)$ with $BRAL = I + T$ and $\mbox{dim} \, \mbox{ker} \, RAL
= \mbox{rank} \, T$. Then $LBRALR = LR + LTR$, whence
\[
(LBR \Pi + I - \Pi) A = I + LTR
\]
(take into account that $A \Pi = \Pi A = A - (I-\Pi)$). This
identity shows that
\[
\mbox{dim} \, \mbox{ker} \, A \le \mbox{rank} \, LTR \le
\mbox{rank} \, T = \mbox{dim} \, \mbox{ker} \, RAL,
\]
whence finally $\mbox{dim} \, \mbox{ker} \, A = \mbox{dim} \,
\mbox{ker} \, RAL$. In the same way one gets $\mbox{dim} \,
\mbox{ker} \, A^* = \mbox{dim} \, \mbox{ker} \, (RAL)^*$. Since
$\mbox{dim} \, \mbox{ker} \, A^* = \mbox{dim} \, \mbox{im} \, A$
for each Fredholm operator $A$, we arrive at assertion $(a)$. \\[1mm]
$(b)$ Let $A_h$ be a limit operator of $A$. Then, by definition,
\[
\|(A_h - V_{-h(n), \, X} A V_{h(n), \, X}) P_{k, \, X}\| \to 0
\quad \mbox{for each} \; k \in \sN.
\]
Thus,
\[
\|R (A_h - V_{-h(n), \, X} A V_{h(n), \, X}) P_{k, \, X} L\| \to 0
\quad \mbox{for each} \; k \in \sN.
\]
Since the projection $\Pi$ commutes with each of the operators
$P_{k, \, X}$, $V_{h(n), \, X}$ and $A$, and since
\[
R V_{h(n), \, X} L = V_{h(n), \, X_N} \quad \mbox{and} \quad R
P_{k, \, X} L = P_{k, \, X_N},
\]
one concludes that
\[
\|(R A_h L  - V_{-h(n), \, X_N} RAL \, V_{h(n), \, X_N}) P_{k, \,
X_N}\| \to 0 \quad \mbox{for each} \; k \in \sN.
\]
Similarly one obtains
\[
\|P_{k, \, X_N} (R A_h L  - V_{-h(n), \, X_N} RAL \, V_{h(n), \,
X_N}) \| \to 0 \quad \mbox{for each} \; k \in \sN.
\]
Thus, $R A_h L$ is the limit operator of $RAL$ with respect to the
sequence $h$. \hfill \qed \\[3mm]
Since the projections $P$ and $\Pi$ also commute, it is an
immediate consequence of the preceding proposition and its proof
that
\[
\mbox{ind}_+ A = \mbox{ind}_+ RAL
\]
and
\[
\mbox{ind}_+ A_h = \mbox{ind}_+ RA_hL = \mbox{ind}_+ (RAL)_h
\]
for each limit operator $A_h \in \sigma_+ (A)$. Thus, the
assertion of Theorem \ref{t10} will follow once we have proved
this theorem for band-dominated operators on $l^p(\sZ, \, X_N)$ in
place of $l^p(\sZ, \, X)$.
\begin{prop} \label{p27}
The assertion of Theorem $\ref{t18}$ holds for all Fredholm
band-domi\-nated operators on $l^p(\sZ, \, X_N)$ (with fixed $N
\in \sN$).
\end{prop}
{\bf Proof.} Let $d < \infty$ be the dimension of $X_N$, and let
$e_1, \, \ldots, \, e_d$ be a basis of $X_N$. Then there are
positive constants $C_1, \, C_2$ such that
\begin{equation} \label{e28}
C_1 \| (x_1, \, \ldots, \, x_d)\|_{l^p} \le \|x_1 e_1 + \ldots +
x_d e_d\|_{X_N} \le C_2 \| (x_1, \, \ldots, \, x_d)\|_{l^p}
\end{equation}
for each vector $(x_1, \, \ldots, \, x_d) \in \sC^d$. Define
$J : l^p(\sZ, \, X_N) \to l^p(\sZ, \, \sC)$ by
\[
(Jx)_{nd+r} := (x_n)_r, \qquad 0 \le r \le d-1
\]
where $(x_n)_r$ refers to the $r$th coordinate of the $n$th entry
$x_n \in X_N$ of the sequence $x$. It follows from (\ref{e28})
that
\[
C_1 \|Jx\|_{l^p(\sZ, \, \sC)} \le \|x\|_{l^p(\sZ, \, X_N)} \le C_2
\|Jx\|_{l^p(\sZ, \, \sC)},
\]
i.e., $J$ is a topological isomorphism from $l^p(\sZ, \, X_N)$
onto $l^p(\sZ, \, \sC)$. The definition of $J$ implies that if $A$
is a Fredholm band operator on $l^p(\sZ, \, X_N)$, then $JAJ^{-1}$
is a Fredholm band operator on $l^p(\sZ, \, \sC)$, and conversely.
Moreover, $\mbox{ind} \, A = \mbox{ind} \, JAJ^{-1}$ in this case.
This identity holds for the plus- and minus-indices as well, since
$J P_{X_N} J^{-1} = P_\sC$. Moreover, one has
\[
J V_{n, \, X_N} J^{-1} = V_{dn, \, \sC} \quad \mbox{and} \quad J
P_{k, \, X_N} J^{-1} = P_{dk, \, \sC}
\]
for all $n \in \sZ$ and $k \in \sN$. These equalities imply that
if $A_h$ is the limit operator of the band-dominated operator $A
\in l^p(\sZ, \, X_N)$ with respect to the sequence $h$, then $J
A_h J^{-1}$ is the limit operator of $JAJ^{-1}$ with respect to
the sequence $dh : \sN \to \sZ, \; m \mapsto d h(m)$, i.e.,
\[
(JAJ^{-1})_{dh} = J A_h J^{-1}.
\]
Summarizing, we obtain
\[
\mbox{ind}_+ \, A = \mbox{ind}_+ \, JAJ^{-1}
\]
and
\[
\mbox{ind}_+ \, A_h = \mbox{ind}_+ \, J A_h J^{-1} = \mbox{ind}_+
\, (JAJ^{-1})_{dh}
\]
for each Fredholm band-dominated operator $A$ on $l^p(\sZ, \,
X_N)$ and for each of its limit operators $A_h \in \sigma_+ (A)$.
Since $dh$ tends to $+ \infty$ whenever $h$ does, one has
$(JAJ^{-1})_{dh} \in \sigma_+ \, (JAJ^{-1})$, and from Theorem
\ref{t10} we infer that $\mbox{ind}_+ \, JAJ^{-1} = \mbox{ind}_+
\, (JAJ^{-1})_{dh}$. Thus, $\mbox{ind}_+ \, A = \mbox{ind}_+ \,
A_h$ for each Fredholm band-dominated operator $A$ on $l^p(\sZ, \,
X_N)$ and for each of its limit operators $A_h \in \sigma_+ (A)$.
The minus-counterpart of this assertion follows analogously. This
proves the proposition and finishes the proof of Theorem
\ref{t18}. \hfill \qed
\section{The Fredholm index of locally compact band-dominated
operators on $L^p(\sR)$} \label{s4}
This section is devoted to the proof of Theorem \ref{t1}. As in
the discrete case, the limit operators approach provides us with a
criterion for the $\hat{\cP}$-Fredholmness of an operator rather
than for its common Fredholmness. Here, $\hat{\cP} =
(\hat{P}_n)_{n \ge 0}$ where $\hat{P}_n : L^p(\sR) \to L^p(\sR)$
is the operator of multiplication by the characteristic function
of the interval $[-n, \, n]$, i.e.,
\[
(\hat{P}_n f)(x) =
\left\{
\begin{array}{ll}
f(x) & \mbox{if} \; x \in [-n, \, n] \\
0    & \mbox{else},
\end{array} \right.
\]
and $\hat{\cP}$-compactness and $\hat{\cP}$-Fredholmness are
defined literally as in the discrete case. The following
proposition can be proved as its discrete counterpart Proposition
\ref{p13}.
\begin{prop} \label{p29}
An operator $A \in L(L^p(\sR))$ of the form $A = I + K$ with
$K \in \cL_p$ is Fredholm if and only if it is $\hat{\cP}$-Fredholm.
\end{prop}
We will now prove Theorem \ref{t1} via a suitable discretization.
Let $\chi_0$ denote the characteristic function of the interval
$[0, \, 1]$. The mapping $G : L^p(\sR) \to l^p(\sZ, \, L^p[0, \, 1])$
which sends the function $f \in  L^p(\sR)$ to the sequence
\[
Gf = ((Gf)_n)_{n \in \sZ} \quad \mbox{where} \quad (Gf)_n :=
\chi_0 U_{-n} f
\]
is a bijective isometry the inverse of which maps the sequence $x
= (x_n)_{n \in \sZ}$ to the function
\[
G^{-1}x = \sum_{n \in \sZ} U_n x_n \chi_0,
\]
the series converging in $L^p(\sR)$. Thus, the mapping
\[
\Gamma : L(L^p(\sR)) \to L(l^p(\sZ, \, L^p[0, \, 1])), \quad
A \mapsto GAG^{-1}
\]
is an isometric algebra isomorphism. It is shown in Proposition 3.1.4
in \cite{RRS4} that
\[
\Gamma (A_h) = (\Gamma (A))_h
\]
for each limit operator $A_h$ of an operator $A \in \cB_p$,
whereas Proposition 3.1.6 in \cite{RRS4} states that $\Gamma$ maps
$\cB_p^\$$ onto $\cA_E^\$$ with $E = l^p(\sZ, \, L^p[0, \, 1])$.
Further, if $A \in L^p(\sR)$ is a locally compact operator, then
the entries of the matrix representation of its discretization
$\Gamma (A)$ are compact operators. Thus, $\Gamma$ maps $\cL_p^\$$
into $I + \cC_E^\$$. Finally, one evidently has
\[
\mbox{ind} \, A = \mbox{ind} \, \Gamma (A)
\]
for each operator $A \in L(L^p(\sR))$, and the Banach space
$L^p[0, \, 1]$ has the $sap$ as already mentioned. Thus, the assertions
of Theorem \ref{t1} follow immediately from their discrete counterparts
Corollary \ref{c15} and Theorem \ref{t18}. \hfill \qed
\section{Applications} \label{s5}
As an application of Theorem \ref{t1}, we are going to examine the
Fredholm properties of operators of the form $I+K$ with $K \in
\cK_p(BUC)$. The latter stands for the smallest closed subalgebra
of $L(L^p(\sR))$ which contains all operators of the form $aCbI$
where $a, \, b \in BUC$ and where $C$ is a Fourier convolution
operator with $L^1$-kernel $k$. Thus,
\[
(Cf)(x) = (k \ast f)(x) = \int_\sR k(x-y) f(y) \, dy, \quad x \in \sR.
\]
In Proposition 3.3.6 in \cite{RRS4} it is verified that
\[
\cK_p(BUC) \subseteq \cL_p^\$.
\]
Hence, Theorem \ref{t1} applies to operators in $\cK_p(BUC)$, and it
yields the following.
\begin{theo} \label{t30}
Let $A \in L(L^p(\sR))$ be a convolution type operator of the form
$I + K$ with $K \in \cK_p(BUC)$. Then \\[1mm]
$(a)$ $A$ is Fredholm if and only if all of its limit operators are
invertible, and if the norms of their inverses are uniformly bounded.
\\[1mm]
$(b)$ if $A$ is Fredholm then, for arbitrary limit operators $B_\pm
\in \sigma_\pm (A)$,
\[
\mbox{\rm ind} \, A = \mbox{\rm ind}_+ B_+ + \mbox{\rm ind}_- B_-.
\]
\end{theo}
One cannot say much about the limit operators of a general
operator $A \in I + \cK_p(BUC)$. It is only clear that they belong
to $I + \cK_p(BUC)$ again. Thus, the computation of the plus- and
minus indices of the limit operators of convolution type operators
will remain a serious problem in general. In what follows we will
discuss some instances where this computation can be easily done
(slowly oscillating coefficients) or is at least manageable
(slowly oscillating plus periodic coefficients).

Let $SO$ stand for the set of all functions $f \in BUC$ which are
{\em slowly oscillating} in the sense that
\[
\lim_{t \to \pm \infty} \sup_{h \in [0, \, 1]} |f(t) - f(t+h)| = 0.
\]
This set forms a $C^*$-subalgebra of $BUC$. Let $\cK_p(SO)$ stand
for the smallest closed subalgebra of $\cK_p(BUC)$ which contains
all operators of the form $aCbI$ where $a, \, b \in SO$ and where
$C$ is a Fourier convolution with $L^1$-kernel. Further, we write
$PER$ for the $C^*$-subalgebra of $BUC$ which consists of all
continuous functions of period 1 on $\sR$. By $\cK_p(PER, \, SO)$
we denote the smallest closed subalgebra of $\cK_p(BUC)$ which
contains all operators of the form $aCbI$ where now $a, \, b \in
PER + SO$ and where $C$ is again a Fourier convolution with
$L^1$-kernel. Similarly, $\cK_p(PER)$ refers to the smallest
closed subalgebra of $\cK_p(BUC)$ which contains all operators
$aCbI$ with $a, \, b \in PER$ and with a Fourier convolution $C$
with $L^1$-kernel.
\begin{lemma} \label{l31}
The limit operators of operators in $\cK_p(SO)$ are operators of
Fourier convolution with $L^1$-kernel, and all limit operators of
operators in $\cK_p(PER, \, SO)$ belong to $\cK_p(PER)$.
\end{lemma}
{\bf Proof.} Operators of convolution are shift invariant with
respect to arbitrary shifts, and operators of multiplications by
functions in $PER$ are invariant with respect to integer shifts.
Hence, operators of this form as well as there sums and products
possess exactly one limit operator, namely the operator itself.
Further, as it has been pointed out in Proposition 3.3.9 in
\cite{RRS4}, all limit operators of operators of multiplication by
slowly oscillating functions are constant multiples of the
identity operator, whence the assertion.
\hfill \qed \\[3mm]
Hence, the determination of the index of a Fredholm operator in $I
+ \cK_p(SO)$ requires the computation of the plus- and the minus
index of an operator of the form $I + C$ where $C$ is a Fourier
convolution with kernel $k \in L^1(\sR)$. Equivalently, one has to
determine the common Fredholm index of operators of the form $I +
\chi_\pm C \chi_\pm I$. The operator $I + \chi_+ C \chi_+ I$ is
the Wiener-Hopf operator with generating function $1+a$ where $a$
is the Fourier transform of the kernel $k$ of $C$. After
reflection at the origin, the operator $I + \chi_- C \chi_- I$
also becomes a Wiener-Hopf operator.

The Fredholm property of Wiener-Hopf operators of this type is
well understood (see \cite{BSi1,GoF1,Kre1}). Since
\[
\lim_{x \to + \infty} a(x) = \lim_{x \to - \infty} a(x) = 0,
\]
one can consider $1+a$ as a continuous function on the one-point
compactification $\dot{\sR}$ of the real line, which is also
called the {\em symbol} of the operator. It turns out that the
Wiener-Hopf operator with symbol $1+a$ is Fredholm if and only if
the function $1+a$ does not vanish on $\dot{\sR}$, and that in
this case its Fredholm index is the negative winding number of the
closed curve $1 + a(\dot{\sR})$ around the origin. This solves the
problem of computing the Fredholm index of an operator in $I +
\cK_p(SO)$ completely and in an easy way.

Let us now turn over to the setting of operators in $I +
\cK_p(PER, \, SO)$. Here we are left with the problem to determine
the Fredholm index of operators of the form $\chi_+ (I + K) \chi_+
I$ on $L^p(\sR^+)$ where $K \in \cK_p(PER)$. The proofs of
Theorems \ref{t1} and \ref{t18} given above offer a way to perform
this calculation. The decisive point is that, due to the
periodicity, the operator
\begin{equation} \label{e32}
\Gamma(\chi_+ (I + K) \chi_+ I) \in L(l^p(\sZ^+, \, L^p[0, \, 1]))
\end{equation}
is a band-dominated {\em Toeplitz operator} the entries of which
are of the form $I + \mbox{\em compact}$ if they are located on
the main diagonal, whereas they are compact when located outside
the main diagonal. Recall that a Toeplitz operator on $l^p(\sZ^+,
\, X)$ is an operator with matrix representation $(A_{i-j})_{i, \,
j \in \sZ^+}$, i.e., the entries of the matrix are constant along
each diagonal which is parallel to the main diagonal.

If now $I + K$ is Fredholm on $L^p(\sR)$, then the Toeplitz
operator (\ref{e32}) is Fredholm, too, and it has the same index.
Employing the reduction procedure used in the proof of Theorem
\ref{t18}, one can further approximate the Toeplitz operator
(\ref{e32}) by a Toeplitz operator on $l^p(\sZ^+, \, \sC^N)$ with
band structure which is also Fredholm and has the same index as
the original operator $I+K$. Thus, we are left with the
determination of the index of a common Toeplitz operator $T(g)$ on
$l^p(\sZ^+, \, \sC^N)$ where each entry $g_{ij}$ of the generating
function $g : \sT \to \sC^{N \times N}$ is a trigonometric
polynomial. This operator can be identified with an operator
matrix $(T(g_{ij}))_{i, \, j = 1}^N$ where each $T(g_{ij})$ is a
Toeplitz band operator on $l^p(\sZ^+, \, \sC) = l^p(\sZ^+)$. As it
is well known (see, e.g., Theorem 6.12 in \cite{BSi1}), this
operator is Fredholm if and only if the common Toeplitz operator
(with scalar-valued polynomial generating function) $T(\mbox{det}
\, g)$ is Fredholm, and the indices of these operators coincide.
Moreover, the index of $T(\mbox{det} \, g)$ is equal to the
negative winding number of the function $\mbox{det} \, g$ with
respect to the origin.

For a general account on matrix functions and the Toeplitz and
Wiener-Hopf operators generated by them, we refer to the
monographs \cite{ClG1} and \cite{LSp1}. Convolution and Wiener-Hopf 
operators with almost periodic matrix-valued generating functions
are thoroughly treated in the monograph \cite{BKS1}. For general 
results about relations between the Fredholmness of a block operator 
and its determinant one should consult Chapter 1 in \cite{Kru1}.

A similar approach is possible for operators in $I +
\cK_p(PER_\sZ, \, SO)$ where $PER_\sZ$ stands for the set of all
functions with integer period. After discretization and
approximation as above, one finally arrives at a block Toeplitz
operator in place of (\ref{e32}) which again can be reduced to a
matrix of Toeplitz operators on $l^p(\sZ^+)$.

The results of Theorems \ref{t1}, \ref{t18} and \ref{t30} can be
completed by an observation made in \cite{Roc9} for the case of
band-dominated operators on $l^p(\sZ, \, \sC)$. This observation
concerns the independence of the Fredholm index on $p$. To make
this statement precise we have to explain what is meant by a
band-dominated operator which acts on different $l^p$-spaces (notice
that the class of all band operators is independent of $p$ whereas
the algebra $\cA_E$ of all band-dominated operators depends on the
parameter $p$ of $E = l^p(\sZ, \, X)$ heavily).

Every infinite matrix $(a_{ij})_{i, \, j \in \sZ}$ induces an
operator $A$ on the Banach space $c_{00}(\sZ, \, X)$ of all functions
$x : \sZ \to X$ with compact support by
\[
i \mapsto (Ax)_i := \sum_{j \in \sZ} a_{ij} x_j.
\]
We say that $A$ {\em extends to a bounded linear operator on}
$l^p(\sZ, \, X)$ or that $A$ {\em acts on} $l^p(\sZ, \, X)$
if $Ax \in l^p(\sZ, \, X)$ for each $x \in c_{00}(\sZ, \, X)$
and if there is a constant $C$ such that $\|Ax\|_p \le C \|x\|_p$
for each $x \in c_{00}(\sZ, \, X)$. If $A$ extends to a
band-dominated operator on both $l^p(\sZ, \, X)$ and
$l^r(\sZ, \, X)$, then we say that $A$ is a {\em
band-dominated operator on} $l^p(\sZ, \, X)$ {\em and}
$l^r(\sZ, \, X)$. Otherwise stated: we consider two band-dominated
operators $B$ and $C$ acting on $l^p(\sZ, \, X)$ and
$l^r(\sZ, \, X)$, respectively, as identical, and we denote them
by the same letter, if their matrix representations coincide.
\begin{prop} \label{p33}
Let $A \in I + \cC_E^\$$ be a Fredholm band-dominated operator both
on $E = l^p(\sZ, \, X)$ and on $E = l^r(\sZ, \, X)$ with $1 < r
< p < \infty$. Then $A$ is a Fredholm band-dominated operator on
each space $l^s(\sZ, \, X)$ with $r < s < p$, and the Fredholm
index $\mbox{\rm ind}_s \, A$ of $A$, considered as an operator
on $l^s(\sZ, \, X)$, is independent of $s \in [r, \, p]$.
\end{prop}
The proof follows exactly the line of the proof of Theorem
\ref{t18}, finally reducing the assertion of the proposition to
the case $X = \sC$ which is treated in \cite{Roc9}. It should be
also mentioned that Proposition \ref{p33} remains valid for
band-dominated operators on $L^p(\sZ^N, \, X)$ with $N$ a positive
integer which also follows from \cite{Roc9}.

In combination with Theorems \ref{t1} and \ref{t30} one gets the
following corollary.
\begin{coro} \label{e34}
$(a)$ Let $A \in I + \cL_q^\$$ be a Fredholm band-dominated operator
both for $q = p$ and for $q = r$ with $1 < r < p < \infty$. Then $A$
is a Fredholm band-dominated operator on each space $L^s(\sR)$ with
$r < s < p$, and the Fredholm index $\mbox{\rm ind}_s \, A$ of $A$,
considered as an operator on $L^s(\sR)$, is independent of $s \in
[r, \, p]$. \\[1mm]
$(b)$  Let $A \in I + \cK_q(BUC)$ be a Fredholm convolution type
operator both for $q = p$ and for $q = r$ with $1 < r < p < \infty$.
Then $A$ is a Fredholm convolution type operator on each space
$L^s(\sR)$ with $r < s < p$, and the Fredholm index $\mbox{\rm ind}_s
\, A$ of $A$, considered as an operator on $L^s(\sR)$, is independent
of $s \in [r, \, p]$.
\end{coro}
{\small Authors' addresses: \\[3mm]
Vladimir S. Rabinovich, Instituto Politechnico National,
ESIME-Zacatenco, Ed.1, 2-do piso, Av.IPN, Mexico, D.F., 07738 \\
E-mail: rabinov@maya.esimez.ipn.mx \\[2mm]
Steffen Roch, Technische Universit\"at Darmstadt,
Fachbereich Mathematik, Schlossgartenstrasse 7, 64289 Darmstadt,
Germany. \\
E-mail: roch@mathematik.tu-darmstadt.de}

\begin{thebibliography}{11}
\bibitem{BKS1}
{\sc A. B\"ottcher, Yu. I. Karlovich, I. M. Spitkovsky}, Convolution
operators and Factorization of Almost Periodic Matrix Functions. --
Birkh\"auser Verlag, Basel 2002.
\bibitem{BSi1}
{\sc A. B\"ottcher, B. Silbermann}, Analysis of Toeplitz operators.
-- Springer-Verlag, Berlin, Heidelberg, New York 1990.
\bibitem{ClG1}
{\sc K. F. Clancey, I. Gohberg}, Factorization of matrix functions
and singular integral operators. -- Birkh\"auser Verlag, Basel 1981.
\bibitem{GoF1}
{\sc I. Gohberg, I. Feldman}, Convolution Equations and Projection
Methods for Their Solution. -- Nauka, Moskva 1971 (Russian, Engl.
transl.: Amer. Math. Soc. Transl. of Math. Monographs, Vol. 41,
Providence, Rhode Island, 1974).
\bibitem{GKr1}
{\sc I. Gohberg, M. Krein}, Systems of integral equations on the
semi-axis with kernels depending on the difference of the
arguments. -- Usp. Mat. Nauk {\bf 13}(1958), 5, 3 -- 72 (Russian).
\bibitem{KaS1}
{\sc N. Karapetiants, S. Samko}, A certain class of convolution
type integral equations  and its applications. -- Izv. Akad. Nauk
SSSR, Ser. Mat. {\bf 35}(1971), 3, 714 -- 726 (Russian).
\bibitem{KaS2}
{\sc N. Karapetiants, S. Samko}, Equations with Involutive
Operators. -- Birkh\"auser Verlag, Boston, Basel, Berlin 2001.
\bibitem{Kre1}
{\sc M. Krein}, Integral equations on the semi-axis with kernels
depending on the difference of the arguments. -- Usp. Mat. Nauk
{\bf 13}(1958), 2, 3 -- 120 (Russian).
\bibitem{Kru1}
{\sc N. Ya. Krupnik}, Banach algebras with symbol and singular
integral operators. -- Shtiintsa, Kishinev 1984 (Russian, English
transl.: Birkh\"auser Verlag, Basel 1987).
\bibitem{LSp1}
{\sc G. S. Litvinchuk, I. M. Spitkovski}, Factorization of measurable
matrix functions. -- Birkh\"auser Verlag, Basel 1987.
\bibitem{RaR5}
{\sc V. S. Rabinovich, S. Roch}, Fredholmness of convolution type
operators. -- {\em In}: Operator Theory: Advances and Applications
{\bf 147}, Birkh\"auser Verlag, Basel, Boston, Berlin 2004, 423 --
455.
\bibitem{RRR1}
{\sc V. S. Rabinovich, S. Roch, J. Roe}, Fredholm indices of
band-dominated operators. -- Integral Equations Oper. Theory
{\bf 49}(2004), 2, 221 -- 238.
\bibitem{RRS1}
{\sc V. S. Rabinovich, S. Roch, B. Silbermann}, Fredholm theory
and finite section method for band-dominated operators. --
Integral Equations Oper. Theory {\bf 30}(1998), 4, 452 -- 495.
\bibitem{RRS2}
{\sc V. S. Rabinovich, S. Roch, B. Silbermann}, Band-dominated
operators with operator-valued coefficients, their Fredholm
properties and finite sections. -- Integral Equations Oper. Theory
{\bf40}(2001), 3, 342 -- 381.
\bibitem{RRS4}
{\sc V. S. Rabinovich, S. Roch, B. Silbermann}, Limit Operators
and Their Applications in Operator Theory. -- Birkh\"auser Verlag,
Basel Boston, Berlin, 2004.
\bibitem{Roc9}
{\sc S. Roch}, Band-dominated operators on $l^p$-spaces: Fredholm
indices and finite sections. -- Acta Sci. Math. (Szeged) {\bf 70}(2004),
3 - 4, 783 -- 797.
\end{thebibliography}
\end{document}